\documentclass[12pt,twoside]{article}

\usepackage{amsmath}
\usepackage{amsfonts}
\usepackage{amssymb}
\usepackage{mathrsfs}
\usepackage{amsfonts}
\usepackage{graphicx}
\usepackage{xcolor}
\usepackage{epstopdf}
\usepackage{amsthm}
\usepackage{cases}
\usepackage{multirow}
\numberwithin{equation}{section}
\usepackage{hyperref}
\makeatletter 
\allowdisplaybreaks
\begin{document}

\date{}
%\author{Qiaohua Liu\thanks{Corresponding author. E-mail: qhliu@shu.edu.cn(Q. Liu),  shendongmei@sfu.edu.cn(D. Shen)}\quad Qian Zhang \quad Dongmei Shen\\
%\small \textit{$^a$Department of Mathematics, Shanghai University, Shanghai 200444, China}}

\author{Qiaohua Liu$^a$\thanks{Corresponding author. E-mail: qhliu@shu.edu.cn(Q. Liu), shendongmei@sfu.edu.cn(D. Shen)}  \; Qian Zhang$^a$\;  Dongmei Shen$^{b*}$\\
\small \textit{$^a$Department of Mathematics, Shanghai University, Shanghai 200444, China}\\
\textit{\small $^b$School of  Statistics and Mathematics, Shanghai Lixin University of Accounting and Finance, }\\ \textit{\small Shanghai, 201209,  China}}

\title{Condition numbers of the mixed least squares-total least squares problem: revisited\thanks{This research is partially supported by the National Natural Science Foundation of China (11001167).}
}
\maketitle
 \hrulefill

\begin{abstract}
A new closed  formula for the first order perturbation estimate   of the mixed least squares-total least squares (MTLS) solution is presented.
It is mathematically equivalent to the one by Zheng and Yang (Numer. Linear Algebra Appl. 2019; 26(4):e2239). With this formula, general and structured normwise, mixed and componentwise condition numbers of the MTLS problem are derived. Perturbation bounds based on the normwise condition number, and
compact forms   for the upper bounds of  mixed and componentwise condition numbers are also given in order for  economic storage and efficient computation.
It is shown that the condition
numbers and perturbation bound of the TLS problem  are unified in the ones of the MTLS problem. \\

\textbf{Keywords}~~ {mixed least squares-total least squares; condition number;  total least squares problem; perturbation bound; linear structure.}
\end{abstract}
\hrulefill
% \linenumbers


\section{ {\textsf{Introduction}}}
\label{intro}
Consider the overdetermined linear system $$Ax\approx b,$$ where $A\in R^{m\times n}, b\in R^m,m\geq n,$ and the matrix $A$
has full column-rank. In some engineering applications,
   the matrix $A$ and the observation $b$  are contaminated
 by some noise, and the total least squares (TLS) model \cite{bj2,gv} is often used to find best approximations to them in Frobenius norm such that
 \begin{equation}
 \min_{E,f} \|[E,f]\|_F \quad
 subject\quad to \quad(A+E)x=b+f.\label{1.1}
  \end{equation}
 The vector $x=x_{\rm TLS}$ satisfying (1.1) is called the TLS solution. If some of rows of $[A,~ b]$ are free of error and  some rows in $[E, f]$ are set zero, the corresponding problem reduces to the total least squares problem with equality constraint \cite{lcz,ljy}.
  If some columns of $A$ are known exactly, and  some columns in $E$, say, the first $n_1$ columns of $E$
are set to be zero,  then the corresponding problem (1.1) is known as the mixed least-squares and total least squares (MTLS) problem \cite{lw,pw,yf}. It arises in the regression analysis \cite{gl}, system identification \cite{ss} and signal processing \cite{vv2}, etc.
 The solution to the MTLS problem is denoted by $x_{\rm M}$. Usually we
 assume that $Ax=b$ is not consistent, otherwise the best minimizer for $[E, f]$ in (\ref{1.1}) is taken to be a zero matrix.\\
\indent For  TLS problems with small  or medium size, a classical direct solver is based on the singular value decomposition (SVD) \cite{bj2,gv,gv2,vv1,xxw2}, where the solution is obtained from the
right singular vector corresponding to the smallest singular value of $[A, b]$.   For the MTLS problem,   QR combined with SVD can be adopted, see e.g.,\cite{vv1}. Recently Liu and Wang \cite{lw} proposed the method
of weighting to interpret the MTLS solution as  a limit of the solution to an unconstrained weighted  TLS (WTLS) problem, as the positive parameter  in the weight matrix tends to zero.   For the solution of problem WTLS, it is    convenient  to apply all   known TLS theories and algorithms  on it. Based on this observation, a QR factorization-based inverse  iteration and
a Rayleigh quotient iteration method were presented in \cite{lw}. The superiority of these methods over the standard QR-SVD algorithm
was demonstrated by numerical experiments.

The condition number is fundamental since it measures the worst-case sensitivity of its solution
to small perturbations in the input data.
The condition numbers of  TLS and its extension the scaled TLS (STLS) problem, by minimizing $\|[E, \lambda f]\|_F$ instead in (\ref{1.1}),
 have been studied widely, e.g. by Zhou et al.\cite{zlw},
  Baboulin and Gratton \cite{bg}, Li and Jia \cite{lj, jl}, Wang et al. \cite{wly}, based on which the explicit expression of the normwise condition number are given.  In \cite{xxw}, Xie et al. provided the perturbation bound for the TLS problem and they also proved that the explicit expression for the absolute normwise condition number in \cite{zlw,bg,lj} are equivalent since they have the same value of the 2-norm.
 Mixed and componentwise condition numbers for a linear function of the solution of the TLS problem were further studied in \cite{ds}. In \cite{zmw,mzw}, Zheng et al. and Meng et al. studied the
condition numbers of multidimensional TLS problems.

Based on the weighting method for problem MTLS \cite{lw}  and the technique \cite{bg} for the conditioning of the standard TLS problem,
 Zheng and Yang \cite{zy}  studied the closed formula for    the first order perturbation estimate of the MTLS solution and gave the
 explicit expressions for
 the condition number of the MTLS problem.  In this paper, we revisit the condition numbers of the MTLS problem.
 We derive another different closed formula  for the first order perturbation estimate of the MTLS solution, and reveal that the formula
 is equivalent to that in \cite{zy}, which also implies the equivalence of the first order perturbation estimates for the TLS problem in \cite{bg,lj}.
General and structured  normwise, mixed and componentwise condition numbers and structured condition numbers  are investigated. We also present the perturbation bound and compact form for the condition numbers to be computable more efficiently with less storage. Numerical tests are given to illustrate our  theoretical results.

 Before our discussion, some notations  are required. $R^m$ and $R^{m\times n}$ denote the spaces of $m\times 1$ and $m\times n$ real matrices, respectively.  $I_n$ denotes the $n
\times n$ identity matrix. $O_{m\times n}$ and  $O_n$ denote   $m\times n$ and $n\times n$ zero matrix, respectively. If subscripts are ignored, the sizes of identity and zero matrices are suitable with context.  $\parallel \cdot \parallel_2$, $\parallel \cdot
\parallel_{\infty}$  and $\parallel \cdot \parallel_{F}$ denote 2-norm, $\infty$-norm  and Frobenius norm of their
arguments, respectively.  For a matrix $A$, $|A|$ is a matrix by taking the absolute value of $a_{ij}$ as elements,  $A^T$ is the transpose of $A$,  $a_j$ is a Matlab notation that denotes
the $j$th column of $A$ and $\sigma_i(A)$ is the $i$-th largest singular value of $A$;
${\rm vec}(A)$ is an operator, which stacks the columns of $A$ one  underneath the other.
  For matrices $A$ and $B$,  the Kronecker product \cite{ls}
of $A$ and $B$ is defined by $A\otimes B=[a_{ij}B]$ and its property is listed as follows \cite{ls,gr}:
$$
\begin{array}l
|A \otimes B |=|A|\otimes |B|,\qquad (A\otimes B)^{T}=A^{T}\otimes B^{T}, \\
{\rm vec}(AXB)=(B^{T}\otimes A){\rm vec}(X).
\end{array}\label{1.01}
$$

\indent Let $x=\phi(a)$ be a  continuous and  Fr\'{e}chet
differentiable mapping from $R^p$ to $R^q$.  For small perturbations $\delta a$,  denote $\delta x=\phi(a+\delta a)-\phi(a)$. According to \cite{ge, gk, rice}, the general normwise condition number
$\kappa(\phi, a)$,  mixed  condition number $m(\phi,a)$ and componentwise condition
number $c(\phi, a)$ are defined and formulated as follows:
$$ \kappa(\phi,a):=\lim_{\varepsilon \rightarrow 0}\sup_{\|\,\delta a\,\|_2\leq\varepsilon\|a\|_2}
\frac{\|\delta x\|_2/\|x\|_2}{\|\delta a\|_2/\|a\|_2}=\frac{\|\phi\prime(a)\|_2\|a\|_2}
{\|\phi(a)\|_2},$$
$$ m(\phi,a):=\lim_{\varepsilon\rightarrow 0}\sup_{\mbox{\tiny $|\delta a|\leq\varepsilon|a|$}}
\frac{\|\delta x\|_{\infty}/\|x\|_{\infty}}{\|\delta a/ a\|_{\infty}}=\frac{\||\phi\prime(a)|\cdot|a|\|_{\infty}}{\|\phi(a)\|_{\infty}},$$
 $$ c(\phi,a):=\lim_{\varepsilon\rightarrow 0}\sup_{\mbox{\tiny $|\delta a|\leq\varepsilon|a|$}}
\frac{\|\delta x/x\|_{\infty}}{\|\delta a/a\|_{\infty}}=\left\|{|\phi\prime(a)|\cdot|a|\over |\phi(a)|} \right\|_{\infty},$$
 where $\phi(a)\not =0$,  $\phi\prime(a)$ denotes the Fr$\acute{e}$chet derivative \cite{ge,rice} of $\phi$ at the point $a$ and $b/a$ is the entry-wise division. Note that $\xi/0$
 is interpreted as zero if $\xi=0$ and infinity otherwise. We assume that  $a$ and $\phi(a)$ have no zero entries throughout this paper.

\section{\textbf{\textsf{A new perturbation estimate of the MTLS solution}}}
\label{sec:1}

For the MTLS problem \cite{vv1}, we assume that the first $n_1$ columns of $A$ are known exactly.  Partition $A = [A_1, A_2]$, where  $A_1\in R^{m \times n_1}$, $A_2 \in R^{m\times n_2}$
and $n_1 + n_2 = n$.  Let the partition $x = [x_1^T, x_2^T]^{T}$  be conformal with the context.  Then the MTLS problem can be formulated as
\begin{eqnarray}
\min_{E_2,f} \|[E_2,f]\|_F \quad subject\quad  to\quad A_1x_1+(A_2+E_2)x_2=b+f.\label{2.1}
\end{eqnarray}
To solve MTLS , a standard way is to factorize $[A, b]$ into the QR form first:
\begin{equation}
Q^{T}[A_1,A_2,b]=\tilde{R}=\begin{array}l\left[\begin{array}{ccc}
 R_{11} & R_{12}&R_{1b}\\
 0 & R_{22} &R_{2b}
 \end{array}
\right] \begin{array}{c}n_1\\ m-n_1\end{array}\\
\quad~ n_1\quad~ n_2\quad~ 1
\end{array}
\label{2.01}
\end{equation}
then solve the reduced TLS problem $R_{22}x_2\approx R_{2b}$ to obtain   $x_2$. The vector $x_1$ is then solved from $R_{11}x_1= R_{1b}- R_{12}x_2$ by backward substitution. According to Golub-Van Loan's theory for standard TLS \cite{gv}, if the genericity condition
\begin{equation}
\sigma_{n_2}(R_{22})> {\sigma}_{n_2+1}([R_{22},R_{2b}])>0,\label{2.20}
\end{equation}
holds,  the reduced TLS problem  and therefore the MTLS problem have a unique solution.

{  {\textbf{Lemma 1}} } \cite{lw}\emph{For the MTLS problem $(2.1)$, let $\tilde{A}=[A,b], W={\rm diag}(O_{n_1},I_{n_2})$, and $\tilde {W}={\rm diag}(W,1)$. Then the MTLS solution $x_{\rm M}$ satisfies the generalized eigenvalue system}
$$
\left[
 \begin{array}{ccc}
 A^{T}A &&A^{T}b  \\
 b^{T}A &&b^{T}b
\end{array}
\right]\left[
 \begin{array}{ccc}
 x_{\rm M}  \\
 -1
\end{array}
\right]=\tilde{\sigma}^{2}_{n_2+1}\left[
 \begin{array}{ccc}
W && \\
  && 1
\end{array}
\right]\left[
 \begin{array}{ccc}
 x_{\rm M}  \\
-1
\end{array}
\right]$$
\textit{where} $\tilde{\sigma}_{n_2+1}^2={\sigma}_{n_2+1}^2([R_{22},R_{2b}])=\lambda_{n+1}(\tilde{A}^{T}\tilde{A}, \tilde{W})$. {\it Here $\lambda_i(M,N)$ denotes the $i$-th generalized
eigenvalue of matrix pair $(M,N)$. }
 \textit{If}
\begin{equation}
\lambda_n(A^TA, W)>\lambda_{n+1}(\tilde{A}^{T}\tilde{A}, \tilde{W}),\label{2.2}
\end{equation}
\textit{then the MTLS  problem} $(2.1)$ \textit{ has a  unique solution  determined by} $x_{\rm M}=( A^{T}A -\tilde{\sigma}^{2}_{n_2+1}W)^{-1}A^{T}b,$
\textit{and it is
equivalent to solving the following optimization problem}
\begin{equation}
\widetilde\sigma_{n_2+1}^2=\min_{x}\frac{\|b-Ax\|_2^2}{1+x^{T}Wx}.\label{2.4}
\end{equation}

 Consider the mapping
  $\varphi:[A,~ b] \longmapsto x_{\rm M}=(A^{T}A-\tilde{\sigma}^{2}_{n_2+1}W)^{-1}A^{T}b$, where $ A$ and $b$ are   input data of the MTLS problem. Let
 $a :={\rm vec}( [A~~ b])$,
then we have the vector representation
$$
x_{\rm M} = \phi(a)=\phi\circ {\rm vec}([A,~ b])=\varphi([A,~ b]).
$$

According to the definitions and formulae of three condition numbers, it is vital to compute the Fr$\acute{e}$chet derivative of $\phi(a)$ for representing the condition numbers of the MTLS problem. To this end, let $\hat {A}=A+\Delta A, \hat{b}=b+\Delta b,$ where $ \Delta A $ and $ \Delta b$
 denote the perturbations to $A$ and $b$ respectively. Consider the perturbed MTLS problem
 \begin{equation}
 \min_{\hat E_2,\hat f}
\|[\hat E_2,\hat f]\|_F \quad subject\quad   to\quad \hat {A}_1 x_1+(\hat{A}_2+\tilde E_2)x_2=\hat{b}+\tilde f. \label{2.5}
\end{equation}

In \cite{zy}, Zheng and Yang proved that
\begin{equation}
\varphi'([A\quad b])=K_{\rm ZY} =[-(x^{T}\otimes D)-(r^T\otimes P^{-1})\Pi_{m,n},\quad D],\label{zy1}
\end{equation}
where  $x$ is the exact MTLS solution, $P=(A^TA-\widetilde\sigma_{n_2+1}^2W),$ $D=P^{-1}(A^T-2{Wxr^T\over \bar\gamma} )$ with $r=Ax-b$, $\bar\gamma=1+x^TWx$ and $\Pi_{m,n}$ is a vec-permutation matrix such that ${\rm vec}(C^T)=\Pi_{m,n}{\rm vec}(C)$ where $C$ is an arbitrary $m\times n$ matrix. Based on this,
they proved that
\begin{equation}
\|K_{\rm ZY}\|_2=\gamma^{1\over 2}\|P^{-1}\Big(A^TA+\gamma^{-1}{\widetilde\sigma_{n_2+1}^2\bar\gamma}(I_n-2{Wxx^TW\over \bar\gamma})\Big)P^{-1}\|_2^{1\over 2},\label{zy2}
\end{equation}
for $\gamma=1+\|x\|_2^2$. Here (\ref{zy2}) is a corrected version of \cite[Eq. (19), Theorem 1]{zy}, where there is a minor error before finishing the final deduction in page 6.

We will adopt a different technique from that in \cite{zy} to derive a new closed formula of  the perturbation estimate.
It is
 clear that when $\|[\Delta A,  \Delta b]\|_{F}$ is sufficiently small, the perturbed MTLS problem has a unique solution $\hat x$ such that $\hat x $ is a real
 analytic function of ${\rm vec}([\Delta A, \Delta b])$ in some neighborhood of the origin.   The following theorem presents a closed formula for the solution of
 the perturbed problem (\ref{2.5}), by   a similar technique to the one  in \cite{lj}.\\

 {\it{\bf Theorem 1}}
\textit{For the MTLS problem (\ref{2.1}) with the genericity condition  (\ref{2.20}) or (\ref{2.2}), denote the  unique solution by $x^{*} = x_{\rm M}$
and define $r=Ax^{*}-b, G(x^*)=[x^{*T},  -1]\otimes I_m$. If $[A,  b]$ is perturbed to $[\hat A,  \hat b]:=[A+\Delta A,  b+\Delta b]$ and $\|[\Delta A,  \Delta b]\|_{F}$ is sufficiently small,
then the perturbed problem (\ref{2.5}) has a unique MTLS solution $\hat x$. Moreover,}
 \begin{equation}
\hat x=x^*+K{\rm vec}([\Delta A, \Delta b]) +{\cal O}(\|[\Delta A, \Delta b]\|_{F}^{2}),\label{2.6}
\end{equation}
where with $P=A^{T}A-\tilde{\sigma}_{n_2+1}^{2}W$, $H_0=I_m- {2rr^T\over \|r\|_2^2}$,
\begin{equation}
K=\varphi'(A,b)=-P^{-1}\left(A^TH_0G(x^{*})+[I_n\otimes r^{T},~ O_{n\times m}]\right).\label{2.7}\\
\end{equation}

{\it Proof.}  For convenience, let $\varepsilon={\rm vec}([\Delta A, \Delta b])$, and $x(\varepsilon)$ be the MTLS solution of the perturbed problem (\ref{2.5}). Similar
to (\ref{2.4}), we can get
$$
x(\varepsilon)={\rm arg}\min_{x}\frac{\|{b+\Delta b}-(A+\Delta A)x\|_2^2}{1+x^{T}Wx}.
$$

It is clear that $x(0)=x^*$, and for sufficiently small $\varepsilon$, $x(\varepsilon)$ is real analytic
 in some neighborhood of the origin. Thus, the Taylor series of $x(\varepsilon)$ with center the origin converges when
 $\varepsilon$ is small enough. So, to prove (\ref{2.6}) it suffices to prove  $\nabla_\varepsilon x(0)$,
 the Jacobian of $ x(\varepsilon)$ at the origin, equals $K$.
To this end, define the two-variable  function$$f(x,\varepsilon)=\frac{\|\hat{b}-\hat{A}x\|_2^2}{1+x^{T}Wx}.$$
The necessary condition for  it to obtain the minimum at $x(\varepsilon)$ is
\begin{equation}
\nabla_xf(x(\varepsilon),\varepsilon)=0.\label{2.71}
\end{equation}
Differentiating  (\ref{2.71}) by $\varepsilon$ with the chain rule gives
$$
\nabla_{x,x}^2f(x(\varepsilon),\varepsilon)\nabla_\varepsilon x(\varepsilon)+\nabla_{\varepsilon,x}^2 f(x(\varepsilon),\varepsilon)=0,
$$
from which it follows that
\begin{equation}
\nabla_\varepsilon x(0)=-(\nabla_{x,x}^2f(x^*,0))^{-1}\nabla_{\varepsilon,x}^2 f(x^*,0),\label{2.8}
\end{equation}
provided that $\nabla_{x,x}^2f(x^*,0)$ is nonsingular.
On the other hand, we note that for $\hat r=\hat Ax-\hat b$,
\begin{equation}
\frac{1}{2}\nabla_xf(x,\varepsilon)=\frac{\hat{r}^{T}\hat{A}}{1+x^{T}Wx}\\
-\frac{\|\hat{r}\|_2^{2}x^{T}W}{(1+x^{T}Wx)^{2}},\label{2.8.1}
\end{equation}
\begin{eqnarray}
\frac{1}{2}\nabla_{x,x}^{2}f(x,\varepsilon)&=&\frac{\hat{A}^{T}\hat{A}}{1+x^{T}Wx}+
4\frac{\|\hat{r}\|_2^{2}Wxx^{T}W}{(1+x^{T}Wx)^{3}} -\frac{\|\hat{r}\|_2^{2}W}{(1+x^{T}Wx)^{2}}
\notag\\&\,& -\frac{1}{(1+x^{T}Wx)^{2}}(2\hat{A}^{T}\hat{r}x^{T}W+
2Wx\hat{r}^{T}\hat{A}),\label{2.80}
\end{eqnarray}
and the fact that $x^{*}=x_{\rm M}$ minimizes  $f(x,0)$,  hence $\nabla_xf(x^{*},0)=0$. Combining this with (\ref{2.8.1}) and (\ref{2.4}), we  get
\begin{equation}
A^{T}r=\frac{\|{r}\|^{2}}{1+x^{*T}Wx^{*}}Wx^{*}=\tilde{\sigma}^{2}_{n_2+1}Wx^{*},\label{2.81}
\end{equation}
 where
$r=Ax^{*}-b$.   Substituting
 $\varepsilon=0$, $x=x^{*}$ and (\ref{2.81}) into (\ref{2.80}), we obtain
\begin{equation}
\frac{1}{2}\nabla_{x,x}^{2}f(x^{*},0)
=\frac{1}{1+x^{*T}Wx^{*}}(A^{T}A-\tilde{\sigma}^{2}_{n_2+1}W).\label{2.8.2}
\end{equation}
 By (\ref{2.2}), we know that $\nabla_{x,x}^{2}f(x^{*},0) $ is positive definite.

 To evaluate $\nabla_{\varepsilon,x}^2 f(x^*,0)$ in (\ref{2.8}), write
  $$
  \hat{r}=\hat{A}x-\hat{b}=\Big([x^{T}, -1]\otimes I_m \Big){\rm vec}([\hat{A},~\hat{b}])=G(x)
{\rm vec}([\hat{A},~\hat{b}])=:G\hat{s}.
$$
Then
 $ \frac{\partial G} {\partial x_i}\hat{s}=\hat{a}_i,i=1,\cdots,n,$
and$$
\frac{1}{2}\nabla_{x}f(x,\varepsilon)=
\frac{1}{1+x^{T}Wx}[\hat{s}^{T}G^{T}\frac{\partial G}{\partial x_1}\hat{s},\cdots,\hat{s}^{T}G^{T}\frac{\partial G}
{\partial x_n}\hat{s}]-\frac{\hat{s}^{T}G^{T}G\hat{s}}{(1+x^{T}Wx)^{2}}x^{T}W,$$
$$\frac{1}{2}\nabla_{\varepsilon, x}^{2}f(x,\varepsilon)=\frac{1}{1+x^{T}Wx}\left[
 \begin{array}{ccc}
 \nabla_\varepsilon(\hat{s}^{T}G^{T}\frac{\partial G} {\partial x_1} \hat{s})\notag\\
 \vdots\\
 \nabla_\varepsilon(\hat{s}^{T}G^{T}\frac{\partial G} {\partial x_n} \hat{s})
\end{array}
\right]-\frac{W}{(1+x^{T}Wx)^{2}}\left[
 \begin{array}{ccc}
 \nabla_\varepsilon(\hat{s}^{T}G^{T}G\hat{s})x_1\notag\\
 \vdots\notag\\
 \nabla_\varepsilon(\hat{s}^{T}G^{T}G \hat{s})x_n\notag
\end{array}
\right],$$
We obtain
$$\frac{1}{2}\nabla_{\varepsilon, x}^{2}f(x,\varepsilon)=\frac{1}{1+x^{T}Wx}\left(
\left[
 \begin{array}{ccc}
 {\hat a_1}^{T}G\notag\\
 \vdots\notag\\
  {\hat a_n}^{T}G\notag
\end{array}
\right]+\left[
 \begin{array}{ccc}
 \hat{r}^{T}\frac{\partial}{\partial x_1}G\notag\\
 \vdots\notag\\
  \hat{r}^{T}\frac{\partial}{\partial x_n}G\notag
\end{array}
\right]\right)-\frac{2Wx\hat{r}^{T}G}{(1+x^{T}Wx)^{2}},$$
in which  $\hat r^T[{\displaystyle \partial G \over\displaystyle \partial x_1}, \cdots, {\displaystyle \partial G \over\displaystyle \partial x_n}]^T=[I_n\otimes \hat r^T, O_{n\times m}]$.
Substituting $\varepsilon =0$ and $x(0)=x_{\rm M}=x^{*} $ into the above equation, by (\ref{2.81}), we get
\begin{eqnarray}
\nabla_\varepsilon x(0)&=&-(\nabla_ {x,x}^{2}f(x^*,0))^{-1}\nabla_{\varepsilon, x}^{2}f(x^*,0)\nonumber\\
&\!=\!&(A^{T}A-\tilde{\sigma}^{2}_{n_2+1}W)^{-1}\left(\frac{2Wx^{*}r^{T}G(x^{*})}{1+x^{*T}Wx^{*}}-A^{T}G(x^{*})-
[I_n\otimes r^{T}\quad O_{n\times m}]\right)\nonumber\\
&\!=\!&(A^{T}A-\tilde{\sigma}^{2}_{n_2+1}W)^{-1}\left(\frac{2A^{T}rr^T}{\|r\|_2^2}G(x^{*})-A^{T}G(x^{*})-
[I_n\otimes r^{T}, O_{n\times m}]\right)\nonumber\\
&\!=\!&(A^{T}A-\tilde{\sigma}^{2}_{n_2+1}W)^{-1}\left(-A^TH_0G(x^{*})-
[I_n\otimes r^{T},  O_{n\times m}]\right).\nonumber
\end{eqnarray}
The proof of the theorem is then completed. \qed\\

{\it  {\bf { Theorem 2}} For the first order perturbation estimate of the MTLS solution $x$, the formulae in (\ref{zy1}) and (\ref{2.7}) are equivalent.}

{\it Proof.} We note that for any $n\times m$ matrix $M_1$,
$$
\begin{array}{rl}
M_1G(x)&=M_1([x^T,\quad -1]\otimes I_m)=M_1[x_1I_m,\cdots, x_nI_m,\quad -I_m]\\
&=[x_1M_1,\cdots,x_nM_1,-M_1]=[x^T\otimes M_1,\quad -M_1],
\end{array}
$$
therefore
$$
\begin{array}{rl}
K&=-P^{-1}(A^TH_0G(x)+[I_n,\quad 0_{n\times 1}]\otimes r^T)\\
&=[-(x^T\otimes (P^{-1}A^TH_0)),\quad (P^{-1}A^TH_0)]-[P^{-1}(I_n\otimes r^T),\quad O_{n\times m}],
\end{array}
$$
where the matrix  $(P^{-1}A^TH_0)$ is exactly the matrix $D$ in  (\ref{zy1}) by  the relation  (\ref{2.81}).

Note that for any $m\times n$ matrix $Y$,
$$
\begin{array}{rl}
 (r^T\otimes P^{-1})\Pi_{m,n}{\rm vec}(Y) &=(r^T\otimes P^{-1}){\rm vec}(Y^T)={\rm vec}( P^{-1}Y^Tr)\\
 &= P^{-1}{\rm vec}(Y^Tr)=P^{-1}{\rm vec}(r^TY)=P^{-1}(I_n\otimes r^T){\rm vec}(Y),
\end{array}
$$
which gives
$$
(r^T\otimes P^{-1})\Pi_{m,n}=P^{-1}(I_n\otimes r^T),
$$
and the assertion in the theorem then follows.\qed\\

\textbf{ {Remark 1}} When $n_1=0$, the MTLS problem (\ref{2.1}) reduces to the standard TLS problem (\ref{1.1})
 and the first order perturbation estimates in (\ref{zy1}) and (\ref{2.6}) become the results in \cite{bg} and  \cite{lj}, respectively, which reveals the equivalence
 of the first order perturbation estimates in \cite{bg} and \cite{lj}.
For $n_2=0$, the MTLS reduces to the
 LS problem, and the estimate in (\ref{zy1}) is exactly the result from \cite{cdw}. Therefore  Theorem 1    unifies the results for the TLS problem and LS problem.

\section{Condition numbers for the MTLS problem}

In this section, we   first consider  condition numbers of the MTLS problem for general matrix $A$ and $b$.
Based on the first order estimate in Theorem 1, we also consider the structured condition numbers   for some MTLS problems with linear structure.

\subsection{Compact  condition numbers and perturbation analysis}

According to   Theorem 1 and the concept of normwise condition number, we   obtain $\phi\prime (a)=K$
and the 2-norm relative condition number of the MTLS problem is given by
\begin{equation}
\kappa(A, b)=\frac{\|K\|_2\|[A, b]\|_{F}}{\|x_{\rm M}\|_2}, \label{3.1}
\end{equation}
where $K$ is defined in (\ref{2.7}).  Note that the expression of $K$ involves Kronecker product, which might  lead to expensive storage and computational cost.
In order to simplify the normwise condition number of  the MTLS problem,
we present the following  theorems.\\

\textit{ {{\bf Theorem 3}} For the MTLS problem (\ref{2.1}), under the genericity condition (\ref{2.20}) or (\ref{2.2}),  the
absolute condition number $\kappa=\|K\|_2$ of the solution $x$ of the MTLS problem has the following   equivalent forms
\begin{eqnarray}
\kappa_1&=&\|P^{-1}(\gamma A^{T}A-A^{T}rx^{T}-xr^{T}A+\|r\|_{2}^{2}I_n)P^{-1}\|_{2}^{\frac{1}
{2}},\label{3.2}\\
\kappa_2&=&\gamma^{1\over 2}\| P^{-1}\Big(A^{T}A+{\gamma^{-1}\widetilde\sigma_{n_2+1}^2\bar\gamma }(I_n-{Wxx^{T}+xx^TW  \over\bar\gamma})\Big)P^{-1}\|_{2}^{\frac{1}
{2}},\label{3.21}
\end{eqnarray}
where $\gamma=1+\|x\|_{2}^{2}$ and $\bar\gamma=1+x^TWx$ with $ W={\rm diag}(O_{n_1},I_{n_2})$. Furthermore,}
\begin{eqnarray}
\kappa_3&=&\|P^{-1}\Big[A^{T},\quad\|x\|_2A^{T}-\|x\|_2\widetilde \sigma_{n_2+1}^2\frac{Wxr^{T}}{\|r\|_{2}^{2}},\quad \|r\|_2I_n-\widetilde\sigma_{n_2+1}^2Wxx^{T}\Big]\|_2,\label{3.33}\\
\kappa_4&=&\left\|P^{-1}\Big[(1+\beta)A^T-\beta\widetilde \sigma_{n_2+1}^2{Wxr^T\over \|r\|_2^2},\qquad \|r\|_2^2I_n-\widetilde \sigma_{n_2+1}^2{Wxr^T}\Big]\right\|_2\label{3.3}
\end{eqnarray}
for $\beta=-1\pm\sqrt{1+\|x\|_2^2}$.

 {\it Proof.}
  For a real matrix $L$, we have $\|L\|_2=\|L^{T}L\|_{2}^{\frac{1}{2}}=\|LL^{T}\|_{2}^{\frac{1}{2}}$. Thus for $K$ we have $$\kappa=\|K\|_{2}=\|KK^{T}\|^{\frac{1}{2}}_2.$$
    Since $P$ is symmetric, by (\ref{2.7}) and with $\bar G=G(x)$, $\Gamma=[I_n, 0_{n\times 1}]\otimes r^{T}$, we observe
    $$
    \bar G\bar G^T=1+\|x\|_2^2=\gamma,\qquad \Gamma\Gamma^T=\|r\|_2^2I_n,\qquad \bar G\Gamma^T=rx^T,
    $$
and the Householder matrix $H_0$ satisfies $H_0r=-r$,   we then  get
\begin{eqnarray}
KK^{T}&=&P^{-1}\left(A^TH_0\bar G+\Gamma\right)\left(\bar G^T H_0^TA+\Gamma^T\right)P^{-1}\notag \\
&=&P^{-1}(\gamma A^{T}A-A^{T}rx^{T}-xr^{T}A+\|r\|_{2}^{2}I_n)P^{-1}.\nonumber
\end{eqnarray}
The relations (\ref{3.2}) and (\ref{3.21}) then follow from  (\ref{2.81}). Furthermore,
\begin{equation}
KK^T=P^{-1}\left([A^{T},\quad I_n]\left[ \begin{array}{ccc}
 \gamma I_m & -rx^{T} \\
 -xr^{T} & \|r\|_{2}^{2}I_n
\end{array}
\right]\left[ \begin{array}{ccc} A\\I_n\end{array}
\right]\right)P^{-1},\label{3.5}
\end{equation}
where with $P_0=I_m-\frac{1}{\|r\|_{2}^{2}}rr^{T}$,
\begin{eqnarray}
&\!\left[\!
 \begin{array}{ccc}
\gamma I_m & -rx^{T} \\
 -xr^{T} & \|r\|_{2}^{2}I_n
\end{array}
\!\right]\!\!=U\left[\!
 \begin{array}{ccc}
I_m+\|x\|_{2}^{2}P_0&O\\
O&\|r\|_{2}^{2}I_n
\end{array}
\!\right]U^T=UD_iD_i^TU^T,\qquad \label{3.6}
\end{eqnarray}
for $i=1,2$ and $D_1={\rm diag}([I_m,~ \|x\|_2P_0], \|r\|_2I_n)$,
$$
D_2=\left[\begin{array}{cc} I_m+\beta P_0& O\\ O& \|r\|_2^2I_n\end{array}\right],\qquad U=\left[\!
 \begin{array}{ccc}
I_m&-\frac{1}{\|r\|_{2}^{2}}rx^{T} \\
O&I_n
\end{array}
\!\right].\label{3.7}
$$
Therefore with $Z_i=[A^T,~ I_n]UD_i$,
$$
\|K\|_2=\|KK^{T}\|_2^{1/2}=\|P^{-1}Z_iZ_i^{T}P^{-1}\|_2^{1/2}=\|P^{-1}Z_i\|_2.\label{3.8}
$$
  By applying the fact in (\ref{2.81}), we obtain the estimates for  $\|P^{-1}Z_1\|_2$ in  (\ref{3.33}) and
$\|P^{-1}Z_2\|_2$ in  (\ref{3.3}). \qed\\

{\it  { {{\bf Theorem 4}}} For the MTLS problem (\ref{2.1}), under the genericity condition (\ref{2.20}) or (\ref{2.2}), if $[A,~ b]$ is perturbed to $[A+\Delta A,~ b+\Delta b]$ and $\|[\Delta A,~ \Delta b]\|_{F}$ is small enough,   then for the exact  MTLS solution $x$ and the solution $\hat x$ to the perturbed problem,  we have for $\Delta x=\hat x-x$ that
\begin{equation}
{\|\Delta x\|_2\over \|x\|_2}\lesssim \kappa_b{\|\Delta b\|_2\over \|b\|_2}+\kappa_A {\|\Delta A\|_2\over \|A\|_2},\label{24}
\end{equation}
where $\kappa_b={\|b\|_2\over \|x\|_2}\Big\|(A^TA-\widetilde\sigma_{n_2+1}^2W)^{-1}A^T\Big\|_2,$ and with $r=Ax-b$,
$$
\kappa_A={\|A\|_2\over \|x\|_2}\Big(\|r\|_2\Big\|(A^TA-\widetilde\sigma_{n_2+1}^2W)^{-1}\Big\|_2+\|x\|_2\Big\|(A^TA-\widetilde\sigma_{n_2+1}^2W)^{-1}A^T\Big\|_2\Big).
$$}

{\it Proof.} Note that
$\Delta x= K{\rm vec}([\Delta A,~ \Delta b])+{\cal O}(\|[\Delta A,~ \Delta b]\|_F^2)$ with $K$ being defined in (\ref{2.7}), where
$$
G(x^*){\rm vec}([\Delta A,~ \Delta b])=\Delta Ax-\Delta b,
$$
$$
[I_n\otimes r^{T},\quad O_{n\times m}]{\rm vec}([\Delta A,~ \Delta b])={\rm vec}(r^T[\Delta A,~\Delta b]\Big[{I_n\atop 0}\Big])=\Delta A^Tr,
$$
therefore
$$
\|\Delta x\|_2\lesssim \|P^{-1}A^T\|_2(\|\Delta A\|_2\|x\|_2+\|\Delta b\|_2)+\|P^{-1}\|_2\|r\|_2\|\Delta A\|_2,
$$
then the result in  (\ref{24}) follows.\\

 {\textbf{\text{Remark 2}}}
   In the case that $n_1=0$, the Kronecker-product-free expression  in (\ref{3.21})  reduces to the  compact formula
  for the normwise
 condition number of the TLS problem  \cite[Eq. (29)]{jl}.
  In   \cite[Theorem 2.3]{wly}, Wang et al. proved that for the TLS problem
$$
\begin{array}{rl}
\bar\kappa_3=\|K\|_2=\|P^{-1}[A^{T},\quad\|x\|_2A^{T}(I_m-\frac{1}{\|r\|_{2}^{2}}rr^{T}),\quad\|r\|_2(I_n-\frac{1}{\|r\|_{2}^{2}}
A^{T}rx^{T})]\|_2,
\end{array}
$$
where $P=A^TA-\sigma_{n+1}([A,~ b])I_n$ and $r=Ax_{\rm TLS}-b$.
Combined with the equality in (\ref{2.81}) for $n_1=0$, the above estimate is just a special case of (\ref{3.33}).

For the perturbation bound of the MTLS problem, the result  in (\ref{24})  is equivalent to Zheng and Yang's result in \cite[Theorem 5]{zy},
and it also    includes  the bound for the standard TLS problem as a special case, see \cite[Theorem 3.1]{xxw}.

 {\textbf{\text{Remark 3}}}  The estimate in (\ref{3.21}) is a little different from that in (\ref{zy2}). We wonder  if   (\ref{3.21}) is equivalent to  (\ref{zy2}) and try to estimate the 2-norm of  $K_{\rm ZY}$ in (\ref{zy1}) with a different approach from that in \cite{zy}:
 $$
 \begin{array}{rl}
 \|K_{\rm ZY}^{\rm new}\|_2^2&=\max\limits_{\|y\|_2=1} \|K_{\rm ZY}^Ty\|_2^2=\max\limits_{\|y\|_2=1}\left\|\left[\begin{array}c -[(x\otimes D^T)-\Pi_{m,n}^T(r\otimes P^{-1})]{\rm vec}(y)\\D^Ty\end{array}\right]\right\|_2^2\\
 &=\max\limits_{\|y\|_2=1}\left\|\left[\begin{array}c -[{\rm vec}(D^Tyx^T)+\Pi_{m,n}^T{\rm vec}(P^{-1}yr^T)]\\D^Ty\end{array}\right]\right\|_2^2\\
 &=\max\limits_{\|y\|_2=1}\Big(\|[\Pi_{m,n}{\rm vec}(D^Tyx^T)+{\rm vec}(P^{-1}yr^T)\|_2^2+\|D^Ty\|_2^2\Big)\\
& =\max\limits_{\|y\|_2=1}\Big(\|{\rm vec}(xy^TD+P^{-1}yr^T)\|_2^2+\|D^Ty\|_2^2\Big)\\
& =\max\limits_{\|y\|_2=1}\Big(\|xy^TD+P^{-1}yr^T\|_F^2+\|D^Ty\|_2^2\Big)\\
& =\max\limits_{\|y\|_2=1}\Big({\rm tr}[(xy^TD+P^{-1}yr^T)(xy^TD+P^{-1}yr^T)^T]+\|D^Ty\|_2^2\Big)\\
& =\max\limits_{\|y\|_2=1}y^T\Big((1+\|x\|_2^2)DD^T+Drx^TP^{-1}+P^{-1}xr^TD+P^{-2}\|r\|_2^2\Big)y,
 \end{array}
 $$
 where  $D=P^{-1}A^TH_0$ by the relation (\ref{2.81}) and
 $$DD^T=P^{-1}A^TAP^{-1},\qquad  Dr=-P^{-1}A^Tr=-\widetilde\sigma_{n_2+1}^2P^{-1}Wx,$$
 from which
 $$
 \|K_{\rm ZY}^{\rm new}\|_2^2=\|P^{-1}\Big(\gamma A^TA-\widetilde\sigma_{n_2+1}^2Wxx^T-\widetilde\sigma_{n_2+1}^2xx^TW+\|r\|_2^2I_n\Big)P^{-1}\|_2,
 $$
where $\gamma=1+\|x\|_2^2$. Consequently  $\|K_{\rm ZY}^{\rm new}\|_2^2$  is equivalent to the estimate  (\ref{3.21}) for  $\|K\|_2^2$. We will provide numerical tests to compare our estimate
with that in \cite{zy} to illustrate our results.

 \textbf{ {\text{Remark 4}}}
In (\ref{2.7}), the matrix $K$  is of size $n\times m(n+1)$,
 while the associated matrices in (\ref{3.21})-(\ref{3.3}) are of size $n\times n$, $n\times (2m+n)$, $n\times (m+n)$ respectively, which is more economic in storage. From the aspect of computation efficiency, the advantages of (\ref{3.21}) over (\ref{3.2}) and (\ref{3.3}) over (\ref{3.33}) are obvious since they require less matrix-product operations.
However, as
pointed out in \cite{bg,hi}, the explicit formulation of matrix cross product $A^{T}A$ and $P^{-1}$ in (\ref{3.21}) is not expected.  The formula in (\ref{3.3}) is preferred in avoiding the
matrix cross prodcut, where in terms of (\ref{2.01}),   its calculation can be implemented by making use of the  intermediate results from solving the MTLS problem. For example,
 the inverse of $P$ can be written as
\begin{equation}
P^{-1}=\left[\begin{array}{cc}
(R_{11}^TR_{11})^{-1}+R_{11}^{-1}R_{12}S^{-1}R_{12}^TR_{11}^{-T}& -R_{11}^{-1}R_{12}S^{-1}\\
-S^{-1}R_{12}^TR_{11}^{-T}&S^{-1}
\end{array}\right],\label{3.16}
\end{equation}
where $S^{-1}=(R_{22}^TR_{22}-\tilde\sigma_{n_2+1}^2I)^{-1}$     can be an intermediate result from solving the TLS problem $R_{22}x_2\approx R_{2b}$ with its normal
equation $x_2=(R_{22}^TR_{22}-\tilde\sigma_{n_2+1}^2I)^{-1}R_{22}^TR_{2b}$, say  in  the  small or medium MTLS, the SVD of $[R_{22}, R_{2b}]$ is available and $S^{-1}$ can be computed
cheaply based on SVD and the result in \cite[Lemma 2]{dwx}; and for solving large MTLS problems when Rayleigh quotient and preconditioned conjugate gradient (RQI-PCG) \cite{bj} method is used, an approximation of
$\widetilde\sigma_{n+1}$ is available and  the linear system $(R_{22}^TR_{22}-\tilde\sigma_{n_2+1}^2I)w=f$   can be efficiently  solved based on preconditioned conjugate gradient method   via  two triangular linear systems in its each iteration step. The computation of  $P^{-1}$ based on triangular linear systems   can therefore be efficiently computed and preserves better numerical stability
  \cite[ch. 8]{hi}.\\

According to the definition of the mixed and componentwise condition numbers, they can be formulated as
$$
m(A,b)={\||K| {\rm vec}([|A|,~|b|])\|_\infty\over \|x\|_\infty},\qquad c(A,b)=\left\|{|K| {\rm vec}([|A|,~ |b|])\over |x|}\right\|_\infty,
$$
where the calculation of $m(A,b), c(A,b)$ involves the Kronecker-product which makes the storage and computation costly. In practical computations,  the upper bounds below are   alternations to improve the computation efficiency. The proof is straightforward.

\label{sec:3}
\textbf{ {Theorem 5}} Under the notation in Theorem 1, the mixed and componentwise condition numbers of the MTLS problem are bounded as
\begin{eqnarray}
m(A,b)&\le& \frac{ \big\|\,|P^{-1}A^TH_0|\Big((|A||x|+|b|)\Big)
+|P^{-1}||A^{T}||r|\big\| _{\infty}}{\|x\|_{\infty}},
\label{3.9}
\end{eqnarray}
\begin{eqnarray}
c(A,b)&\leq& \Bigg\|\frac{|P^{-1}A^TH_0|(|A||x|+|b|)
+|P^{-1}||A^{T}||r|}{|x|} \Bigg\|_\infty.
\label{3.10}
\end{eqnarray}

\subsection{Structured condition numbers }

 If the matrix $A$ lies in a linear subspace ${\cal S}$ which consists of a class of structured matrices, then  any matrices in  ${\cal S}$ can be represented by a linear combination of
 a linearly independent matrices $S_1, S_2, \cdots, S_q\in {\cal S}$, i.e.
$
 A=\sum\limits_{i=1}^q\alpha_iS_i$. Then
 ${\rm vec}(A)=\sum\limits_{i=1}^q \alpha_i{\rm vec}(S_i)=\Phi_A^{\rm struct}\alpha$, where $\Phi_A^{\rm struct}=[{\rm vec}(S_1),{\rm vec}(S_2),\cdots, {\rm vec}(S_q)]$ and
$\alpha=[\alpha_1,\alpha_2,\cdots, \alpha_q]^T$ and by the statement in \cite[Theorem 4.1]{lj},  $\Phi_{A}^{\rm struct}$  is column orthogonal and has  full column rank,  with at most one
nonzero entry in each row.

 Note that
 $$
 {\rm vec}([A, b])=\Phi_{A, b}^{\rm struct}s:=\left[
 \begin{array}{cc}
 \Phi_{A}^{\rm struct}&0\\
 0& I_m\end{array}\right]\left[\begin{array}c \alpha\\  b\end{array}\right],
 $$
and for the perturbed MTLS problem, if we restrict the perturbation matrices $[\Delta A, \Delta b]$ to have
the same structure as that of $[A, b]$, that is, ${\rm vec}([\Delta A, \Delta b])=\Phi_{A, b}^{\rm struct}\epsilon$ where $\epsilon\in R^{q+m}$.

Define the mapping $\phi$ from $R^{q+m}$ to $R^{n}$ such that
$
\phi(\big[{\alpha\atop b}\big])=x_{\rm M}=(A^TA-\widetilde\sigma_{n_2+1}^2W)^{-1}A^Tb.
$
Based on (\ref{2.6}), the first order perturbation result becomes
$\Delta x=K\Phi_{A, b}^{\rm struct}\epsilon+{\cal O}(\|\epsilon\|_2^2).$
According to the concept of condition numbers, the relative norwise, mixed and componentwise condition numbers for structured MTLS take  following forms
\begin{equation}
\begin{array}l
\kappa^{\rm struct}(\alpha, b)=\|K\Phi_{A, b}^{\rm struct}\|_2{\displaystyle\|[\alpha^T, b^T]\|_2\over\displaystyle \|x_{\rm M}\|_2},\\[10pt]
m^{\rm struct}(\alpha, b)={\displaystyle\left\| |K\Phi_{A, b}^{\rm struct}|\cdot [|\alpha|^T, |b|^T]^T\right\|_\infty\over\displaystyle \|x_{\rm M}\|_\infty},\\[10pt]
c^{\rm struct}(\alpha, b)=\left\|{\displaystyle |K\Phi_{A, b}^{\rm struct}|\cdot [|\alpha|^T, |b|^T]^T|\over\displaystyle |x_{\rm M}|}\right\|_\infty,
\end{array}\label{3.11}
\end{equation}
where
\begin{equation}
\begin{array}{rl}
K\Phi_{A, b}^{\rm struct}&=-P^{-1}\left(A^TH_0[x^T\otimes I_m, -I_m]+[I_n\otimes r^{T},~ O_{n\times m}]\right)\left[\begin{array}{cc}\Phi_{A}^{\rm struct} &0\\0 & I_m\end{array}\right]\\
&=-P^{-1}\Big(A^TH_0[S_1x,\cdots, S_qx, -I_m]+[S_1^Tr,\cdots, S_q^Tr, O_{n\times m}]\Big),
\end{array}\label{3.12}
\end{equation}
which is Kronecker product-free and can be computable
more efficiently with less storage.

\section{\textbf{{Numerical experiments}}}
\label{sec:4}
 In this part, we first present numerical experiments to verify the  utility of our first order perturbation results, and then compare three types of condition numbers from tests. All experiments are coded by MATLAB R2012b  with machine precision $2.22\times10^{-16}$ and done on a PC with Intel Core i5-5200U CPU @ 2.20GHz  and
the memory is 4 GB.

{\bf Example 1} Consider the estimation of the parameters in a transfer function model \cite{you}, given in its errors-in-variables form
$$
\begin{array}l
C(q^{-1})y_0(t)=B(q^{-1})u_0(t),\\
u(t)=u_0(t)+\Delta u(t),\\
y(t)=y_0(t)+\Delta y(t),
\end{array}
$$
where $u_0(t)$ and $y_0(t)$ are the unmeasurable noise-free inputs and outputs, $u(t)$ and $y(t)$ the noisy measured inputs and outputs, while $\Delta u(t)$ and $\Delta y(t)$ represent all stochastic disturbances to the inputs and outputs, respectively; $A(q^{-1})$ and $B(q^{-1})$ are polynomials taking the form
$$
\begin{array}l
B(q^{-1})=b_1q^{-1}+\cdots+b_{n_1}q^{-n_1},\\
C(q^{-1})=1+c_1q^{-1}+\cdots+c_{n_2}q^{-n_2},
\end{array}
$$
and $q^{-1}$ is a backward shift operator such that $q^{-1}y(t)=y(t-1)$. In order to estimate the parameters in transfer function,
we need to solve the approximate system $\phi(t)^Tx\approx y(t)$, for $t=1,2,\cdots, N$ and
$$
\phi(t)=[u(t-1), \cdots, u(t-n_1), -y(t-1),\cdots, -y(t-n_2)]^T,   \qquad x=[b_1,\cdots, b_{n_1}, c_1,\cdots, c_{n_2}]^T,
$$
in which the entries $u(j), y(j)$ with $j\le 0$ are set to be zero. This leads to the set of linear equations $Ax\approx z$ for
$$
A=[\phi(k_0+1),\cdots, \phi(k_0+m)]^T, \qquad z=[y(k_0+1),\cdots, y(k_0+m)]^T.
$$
with  $k_0+1(\ge 1)$ and $m$ being the starting point  and   the number of chosen samples for parameter estimation.
In \cite{vv2},  Van Huffel  and  Vandewalle proposed the TLS model to solve the linear equation when the inputs $u(t)$ and outputs $y(t)$ are affected by the noise. When the inputs $u(t)$ are noise-free, i.e. $\Delta u(t)=0$, the MTLS model should be used to solve the linear equation. In our test, we assume that the inputs $u(t)$ are noise-free and the outputs $y(t)$ are affected by white noise with zero mean and
variance $0.01$.

\renewcommand\tabcolsep{10.0pt}
\begin{table}
\caption{Comparisons of   the absolute normwise condition number with  different forms for perturbed MTLS problems}
\label{tab1}       % Give a unique label
\begin{tabular}{ccccccccc}
\\\hline
\multirow{2}{*}{$\epsilon$}&\multirow{2}{*}{$\|\Delta x\|_2$}& \multirow{2}{*}{$\eta_{\Delta x}^{\rm ZY}$}&\multirow{2}{*}{$\eta_{\Delta x}^{\rm new}$}&
 $\|K_{\rm ZY}^0\|_2$&$\|K_{\rm ZY}\|_2$& $\|K^0\|_2$& $\kappa_2$&$\kappa_4$\\
 &&&&(\ref{zy1})&(\ref{zy2})&(\ref{2.7})& {(\ref{3.21})}& {(\ref{3.3})}\\ \hline
 1e-2&   3.63e-2& 2.93e-4&2.93e-4&18.50&19.80&18.50&18.50&18.50\\
 1e-4&    1.77e-3&2.00e-7&2.00e-7&11.43&16.46&11.43&11.43&11.43\\
 1e-6&  1.95e-5&7.89e-11&7.89e-11&50.51&56.79&50.51&50.51&50.51\\
 1e-8&    2.20e-8&5.41e-15&5.41e-15&13.56&14.13&13.56&13.56&13.56\\\hline
\end{tabular}
\end{table}

Take $m=30, n=20, n_1=10$ and $k_0=\max(n_1, n_2)=10$ and  generate entrywise perturbation
$$
[\Delta A,~ \Delta b]=\epsilon\cdot\textsf{rand}(m,n+1)\odot [A, ~b],
$$
where $\odot$ denotes the entrywise multiplication.
Denote
$$
\eta_{\Delta x}^{\rm ZY}={\|\Delta x-K_{ZY}^0{\rm vec}([\Delta A,~ \Delta b])\|_2},\quad
\eta_{\Delta x}^{\rm new}={\|\Delta x-K^0{\rm vec}([\Delta A,~ \Delta b])\|_2},
$$
where $K_{\rm ZY}^0, K^0$ are computed from   formulae  (\ref{zy1}) and (\ref{2.7}) respectively.
For different $\epsilon $ and random perturbations, in Table 1
we compare the first order perturbation estimates in (\ref{zy1}), (\ref{2.7}) of the MTLS solution and  the estimates in  (\ref{zy2}),(\ref{3.21}) and (\ref{3.3}) for  the absolute normwise condition number.

In Table \ref{tab1}, we observe that for different parameter $\epsilon$ and perturbations, $\eta_{\Delta x}^{\rm ZY}, \eta_{\Delta x}^{\rm new}$ are all of the magnitude ${\cal O}(\epsilon^2)$, indicating the first order perturbation estimates in (\ref{zy1}) and (\ref{2.7}) are both correct. Among five methods for evaluating the normwise condition number
$\|K\|_2$ or $\|K_{\rm ZY}\|_2$, we find that four methods give the same value while the estimate via (\ref{zy2}) does not match the true value of $\|K_{\rm ZY}\|_2$ in (\ref{zy1}),
which illustrates our theoretical results.\\

{}

\textbf{Example 2} In this example, we compare the relative error of the MTLS solution with the estimated upper bounds  based on the normwise condition number or the perturbation bound in (\ref{24}) as well. Firstly, we construct the random MTLS problems as follows. The coefficient matrices $A$ and $b$ are generated according to the QR factorization (\ref{2.01}),
where $Q$ is a random orthogonal matrix, and
  $$
  [R_{11},R_{12},R_{1b}]=\texttt{triu}(\texttt{qr}(\texttt{rand}(n_{1},n+1))).
  $$
Here \texttt{qr}(..), \texttt{triu}(..) are Matlab commands to produce QR factorization and the upper triangular part of a matrix, respectively.
With random unit vectors $y \in R^{m}$ and $z \in R^{n+1}$, set $Y=I_m-2yy^{T},Z=I_{n+1}-2zz^{T},$ $D={\rm diag}(n_2,n_2-1,\cdots,1,1-e_p)$ for given parameter $e_p$, the subblock matrix  $[R_{22},
R_{2b}]$ is generated by
$$
[R_{22},R_{2b}]=\texttt{triu}\left(\texttt{qr}\left(Y\left[
 \begin{array}{ccc}
D\\O
\end{array}
\right]Z^{T}\right)\right),
$$
which is similar to that in \cite{bg}.  Due to the interlacing property \cite{bj2}, we get
$$
0\le \sigma_{n_2}(R_{22})-\sigma_{n_2+1}([R_{22}, R_{2b}])\leq \sigma_{n_2}([R_{22}, R_{2b}])-
\sigma_{n_2+1}([R_{22}, R_{2b}])=e_p,
$$
where the quantity $\sigma_{n_2}(R_{22})-\sigma_{n_2+1}([R_{22}, R_{2b}])$ measures the distance of our problem to nongenericity.
By varying $e_p$, $m$, $n$ and $n_1$, we can generate different MTLS
 problems.  With small values of  $e_p$, it is possible to study the utility of  condition numbers of new form.

Consider the perturbation as in Example 1 with $\epsilon=10^{-10}$. In Table 2, we compare the efficiency of computing the normwise condition number with different forms for the MTLS problem.
 Let $\kappa_0=\|K^0\|_2$   be computed via formula (\ref{2.7}) and  denote  the parameter $\varepsilon_1$ and relative normwise condition numbers of the MTLS problem as
 \begin{equation}
\varepsilon_1={\|[\Delta A,~ \Delta b]\|_F\over \|[A,~ b]\|_F},\qquad \tilde\kappa_{i}=\frac{\kappa_{i}\|\,[A,~ b]\,\|_F}{\|x\|_2}, \qquad i=0, 4,
\end{equation}
where during the computation of $\kappa_i$, the matrix inverse $P^{-1}$ is computed  based on the formula   (\ref{3.16}) and  \cite[Lemma 2]{cdw}.

\renewcommand\tabcolsep{12.0pt}
\begin{table}
\caption{Comparisons of the upper  bounds with forward errors for  perturbed MTLS problems}
\label{tab2}       % Give a unique label
\begin{tabular}{cccccccc}
\hline\noalign{\smallskip}
{$m=300$}&{$n=200$}&{${\|\Delta x\|_2\over \|x\|_2}$ }&$\varepsilon_1\tilde\kappa_0$(Time) &$\varepsilon_1\tilde \kappa_4$(Time)  & (\ref{24})(Time)\\
\noalign{\smallskip}\hline\noalign{\smallskip}
$e_p=0.9$&$n_1$=60&  1.14e-9&3.90e-7(14.23)&3.90e-7(7.14)&   4.45e-7(7.14)\\
$\quad$&$n_1$=120& 1.28e-9&1.41e-6(15.32)&1.41e-6(7.69)&   1.69e-6(7.69)\\
$\quad$&$n_1$=180& 3.89e-9&   1.11e-5(14.01)& 1.11e-5(7.04)&1.28e-5(7.04)\\
\noalign{\smallskip}\hline\noalign{\smallskip}
$e_p$=0.0009&$n_1$=60& 9.26e-8&   5.89e-5(14.46)&   5.89e-5(7.29)&   8.98e-5(7.29)\\
$\quad$&$n_1$=120& 3.19e-7&1.52e-4(13.92)&1.52e-4(6.99)&1.64e-4(7.02)\\
$\quad$&$n_1$=180&  4.92e-8&   3.13e-5(14.01)&   3.13e-5(7.05)&3.27e-5(7.05)\\
\noalign{\smallskip}\hline\noalign{\smallskip}
\end{tabular}
\end{table}

In Table \ref{tab2}, we listed the upper bounds for the relative forward error of the MTLS solution, where   the CPU time for seconds to compute   corresponding upper bounds is deplayed.
  It is observed from Table 2 that the computation of $\tilde\kappa_0$ is much less efficient due to the large size of $K$,  and the bound in (\ref{24}) and upper bounds based on
   $\tilde \kappa_4$ are very tight and  efficient.\\

 {\bf Example 3}  Consider the intercept model arising in the ``errors in variables'' regression model \cite{gl}:
 $ \alpha+x_1c_1+\cdots +x_nc_n=b,$
where $b$ and $c_i$ are observed $m\times 1$ vectors, $\alpha\in R^m$ is the intercept vector of the linear model. The model gives rise to the overdetermined  set of equations
 $$
 [1_m, C]\left[\begin{array}c\alpha\\ x\end{array}\right]=b,
$$
for $1_m=[1,\cdots, 1]^T$ and $C=[c_1,c_2,\cdots, c_n]\in R^{m\times n}$, in which the first column of the left-hand side matrix
is known exactly,  leading to the mixed LS-TLS problem.

In this example, we choose $b$ as a random vector, while for the matrix $C$, we consider the following two cases:

(a) Take $m=6$, $n=4$  and
$$
C=\left[
 \begin{array}{cccc}
 \delta & 0&0&0  \\
 0&\delta& 0&0\\
 0&0&0&1\\
 0&0&1&0\\
 0&0&1&0\\
 0&0&0&\delta
\end{array}
\right], $$
 where $\delta$ is a tiny positive parameter.

Set $A=[1_m, C]$ and in the perturbed MTLS problem, the first column of $A$ is not perturbed, and the matrix $C$ is perturbed according to its structure, i.e.,
$$
[ \Delta C,~~ \Delta b] =10^{-10}\cdot  \textsf{rand}(m,n+1)\odot [C,~~ b].
 $$

{\iffalse

\textbf{Example 4}
  In this example, we construct the matrix $A$ and  vector $b$  as follows,
$$
A=\left[
 \begin{array}{cccc}
 \delta & 0&0&0  \\
 0&\delta& 0&0\\
 0&0&0&1\\
 0&0&1&0\\
 0&0&1&0\\
 0&0&0&\delta
\end{array}
\right] \in R^{6\times 4}$$
and $b \in R^{6\times 1}$ is a random vector,  where $\delta$ is a tiny positive parameter.
Consider the perturbation
$$
 \Delta A =\epsilon\cdot\Delta A_{1}\odot A,  \, \Delta b=\epsilon\cdot\Delta b_1\odot b,
 $$
 where $\epsilon=10^{-10}$, $\odot$ denotes the Hadamard product and each components of $\Delta A_{1} \in R^{6\times 4 }$ and $\Delta b_1 \in R^{6}$ are uniformly distributed in the
 interval $(0,1)$. \fi}

 Define
 $$\varepsilon_2=\min \{\epsilon: |\Delta A|\le \epsilon |A|,  |\Delta b|\le \epsilon |b|\}.
 $$
In Table 3,  we compare the relative perturbations $\frac{\|\Delta x\|_2}{\|x \|_2}$ ,$\frac{\|\Delta x\|_{\infty}}{\|x \|_{\infty}}$ and $\|\frac{\Delta x}{x}\|_{\infty}$ with the backward error multiplied by our normwise condition number, the upper bounds in (\ref{3.9})-(\ref{3.10}) (denoted by $m^u, c^u$, resp.) of mixed and componentwise condition numbers (denoted by $m, c$, resp.). The tabulated results show
that the normwise condition number multiplied by backward error is far from the true value of relative error of the solution when $\delta$ decreases, while the mixed and componentwise condition numbers based bounds can
estimate the forward error much more tightly.
Moreover, the upper bounds of mixed and componentwise condition numbers are very sharp but with less matrix operations.

\renewcommand\tabcolsep{6.0pt}

 {\small
 \begin{table}
\caption{Comparision of upper bounds with  forward errors  for perturbed MTLS problems}
\label{tab:3}       % Give a unique label

\begin{tabular}{llllllllll}
\hline\noalign{\smallskip}
$\delta$&$n_1$&$\frac{\|\Delta x\|_2}{\|x\|_2}$&$\varepsilon_1 \tilde\kappa_4$&$\frac{\|\Delta x\|_{\infty}}{\|x\|_{\infty}}$&$\varepsilon_2 m$&$\varepsilon_2 m^u$&$\|\frac{\Delta x}{x }\|_{\infty}$&$\varepsilon_2 c$&$\varepsilon_2 c^u$\\
\noalign{\smallskip}\hline\noalign{\smallskip}
$10^{-2}$&1&  9.68e-11&1.52e-8&9.34e-11&3.64e-9&3.89e-9&3.32e-10&1.25e-8&1.38e-8\\
$\quad$&3& 5.59e-11&2.63e-8&5.59e-11&2.61e-10&2.82e-10&7.41e-11&1.88e-9&2.50e-9\\
\noalign{\smallskip}\hline\noalign{\smallskip}
$10^{-4}$&1& 2.09e-10&7.55e-6&2.06e-10&1.85e-9&4.47e-9&2.96e-10&2.65e-9&6.23e-9\\
$\quad$&3&  3.42e-11&2.85e-6&2.98e-11&4.58e-10&5.33e-10&9.27e-11&1.42e-9&1.66e-9\\
\noalign{\smallskip}\hline\noalign{\smallskip}
 $10^{-6}$&1&  1.09e-10&2.73e-4&9.93e-11&9.01e-10&9.75e-10&4.02e-10&2.79e-9&3.08e-9\\
$\quad$&3&   1.00e-11&5.91e-4&7.69e-12&1.06e-9&1.24e-9&2.14e-10&1.81e-8&2.12e-8\\
\hline\noalign{\smallskip}
\end{tabular}
\end{table}}

(b) Take $C=T$ to be a large $m\times (m-2\omega)$ Toeplitz matrix, whose first column is given by
$
%t_{i,1}={1\over \sqrt{2\pi\tau^2}}{\rm exp}\Big[{-(\omega-i+1)^2\over {2\pi\tau^2}}\Big],\quad i=1,2,\cdots, 2\omega+1,
t_{i,1}=i\ \mbox{for}\ i=1,2,\cdots, 2\omega+1,
$
and zero otherwise. Entries in the first row are given by $t_{1,j}=t_{1,1}$ if $j=1$, and zero otherwise.
% where $\tau=1.25, \omega=8$.

For the coefficient matrix $A=[1_m, T]$ in the intercept model, it has the special structure such that
$$
A=S_1+t_{1,1}S_2+\cdots+t_{2\omega+1,1}S_{2\omega+2},
$$
where $S_1=[1_m, O_{m\times(m-2\omega)}], S_i=[0_{m\times 1}, \hat S_{i-1}]$ for $i=2,\cdots, 2\omega+2$  and  $\hat S_{1}=\Big[{I_{m-2\omega}\atop O_{2\omega\times (m-2\omega)}}\Big]$
$\hat S_i=Y_0\hat S_{i-1}$  and $Y_0$ is a  lower shift matrix of order $m$. Note that in (\ref{3.12}), for $q=2\omega+2$,
$$
\begin{array}l
[S_1x, S_2x,\cdots, S_{2\omega+2}x]=[x_11_m, T_x],\\[10pt]
[S_1^Tr, S_2^Tr,\cdots, S_{2\omega+2}^Tr]={\rm diag}((1_{m}^Tr), H_r),
\end{array}
$$
where $T_x$ is an $m\times (m-2\omega)$ lower Toeplitz  matrix with $x_2,\cdots, x_{m-2\omega+1}, 0,\cdots, 0$ in its first column, and $H_r$ is an  $(m-2\omega)\times (2\omega+1)$ anti-upper Hankel matrix with $r_1,\cdots, r_{m-2\omega}$ in its first column and $r_{m-2\omega},\cdots, r_m$ in its last row.

Let $b=\lambda \bar b$ with $\bar b$ being a random vector, consider the random entrywise perturbation
 $$
\Delta t_1 = 10^{-10}\cdot {\sf randn}(m, 1)\odot t_{1},\quad \Delta b=10^{-10}\cdot {\sf randn}(m,1)\odot b,
 $$
 such that $T$ is perturbed to  a matrix with the same structure as $T$.

To evaluate  sharp bounds via structured condition numbers (denoted by $\kappa^{s}, m^s, c^s$, resp.) in (\ref{3.11})-(\ref{3.12}), denote
$$
\varepsilon_1^{s}={\|[\Delta \alpha^T, \Delta b^T]\|_2\over \|[\alpha^T, b^T]\|_2},\quad
\varepsilon_2^{s}=\min\{\epsilon: |\Delta \alpha|\le \epsilon |\alpha|, |\Delta b|\le \epsilon |b|\},
$$
in which $\alpha=[1, t_{1,1},\cdots, t_{2\omega+1,1}]^T$.  In the experiment, we take $n_1=1, \omega=8$ and $m=500, 1000$ to compare   general condition numbers with structured ones in Table \ref{tab:4}, where for general mixed and componentwise condition numbers, we use $m^u, c^u$ instead to evaluate the upper bounds of ${\|\Delta x\|_\infty\over \|x\|_\infty}$ or $\|{\Delta x\over x}\|_\infty$, since for large values of $m$, the storage and computation of matrix $K$ in (\ref{2.7}) might exceed the  maximal  memory of the computer.

The tabulated results in Table \ref{tab:4} show that  for evaluating the forward error ${\|\Delta x\|_2\over \|x\|_2}$ of large MTLS problems,
 the estimates  based on structured normwise condition numbers are tighter than those based on unstructured ones,  but both become worse when the matrix $[A,  b]$ becomes badly scaled. The estimates based on general and structured mixed or componentwise condition number are sharp to evaluate ${\|\Delta x\|_\infty\over \|x\|_\infty}$ or $\|{\Delta x\over x}\|_\infty$, they are only one or two  orders of magnitude higher than corresponding relative forward errors, indicating that
mixed and componentwise condition numbers are more robust and preferred  for badly scaled problems.

\renewcommand\tabcolsep{8.0pt}
\begin{table}
\caption{Comparision of upper bounds with  forward errors  for perturbed MTLS problems}
\label{tab:4}
\begin{tabular}{llllllllll}
\hline\noalign{\smallskip}
\multicolumn{9}{l}{$m=500$}\\\hline
$\lambda$&$\frac{\|\Delta x\|_2}{\|x\|_2}$&$\varepsilon_1 \tilde\kappa_4$&$\varepsilon_1^{ s}\kappa^{  s}$&$\frac{\|\Delta x\|_{\infty}}{\|x\|_{\infty}}$&$\varepsilon_2 m^u$&$\varepsilon_2^{  s}m^{  s}$&$\|\frac{\Delta x}{x }\|_{\infty}$&$\varepsilon_2 c^u$&$\varepsilon_2^{ s}c^{ s}$\\\hline
$10^{-2}$&  4.6e-9&9.6e-7&4.3e-8&4.4e-9&1.0e-7&2.5e-8&4.4e-6&1.1e-4&2.8e-5\\
$10^{-4}$& 5.2e-9&8.9e-5&4.1e-6&6.2e-9&8.7e-8&2.7e-8 &1.5e-6&2.4e-5&1.2e-5\\
$10^{-6}$&   5.7e-9&8.6e-3&3.9e-4&5.6e-9&1.1e-7&2.7e-8&4.5e-7&1.1e-5&3.8e-6\\
\noalign{\smallskip}\hline\noalign{\smallskip}
\multicolumn{9}{l}{$m=1000$}\\\hline
$\lambda$&$\frac{\|\Delta x\|_2}{\|x\|_2}$&$\varepsilon_1 \tilde\kappa_4$&$\varepsilon_1^{  s}\kappa^{  s}$&$\frac{\|\Delta x\|_{\infty}}{\|x\|_{\infty}}$&$\varepsilon_2 m^u$&$\varepsilon_2^{  s}m^{ s}$&$\|\frac{\Delta x}{x }\|_{\infty}$&$\varepsilon_2 c^u$&$\varepsilon_2^{  s}c^{  s}$\\\hline
$10^{-2}$& 5.2e-9&  1.4e-6&4.2e-8&4.0e-9&1.9e-7&3.9e-8&1.9e-6&6.4e-5&2.7e-5\\
$10^{-4}$& 1.2e-8& 1.0e-4&3.3e-6&1.2e-8&2.4e-7&5.9e-8&4.6e-6&1.7e-4&6.0e-5\\
$10^{-6}$& 1.5e-8& 1.3e-2&4.2e-4&9.7e-9&1.5e-7&3.6e-8&8.4e-6&1.5e-4&4.3e-5\\
\noalign{\smallskip}\hline\noalign{\smallskip}
\end{tabular}
 \end{table}

\section{Conclusion}
\indent In this paper, we present a new closed formula for the first order perturbation estimate of the MTLS solution, which is also illustrated to be equivalent to that in \cite{zy}. Normwise, mixed/componentwise condition numbers and corresponding structured condition numbers of the MTLS problem are also derived.  These expressions all involve matrix Kronecker product operations, we  propose different skills to simplify the  expressions to improve the computational efficiency.  For the normwise condition number, we show that it can be transformed into several compact forms, and the perturbation bound is also an alternation.
From a number of numerical tests, we can see that the new forms and bounds for the normwise condition number   have great computational efficiency and require less storage.
For structured/sparse and badly-scaled  MTLS problems,  the  sparse pattern and the magnitude of
the entries are  better utilized in (structured) mixed and componentwise condition numbers than normwise condition number,  and it is more suitable to adopt mixed/componentwise condition numbers to measure the conditioning of  badly scaled or structured MTLS problems.

%\bibliography{mybibfile}

\end{document}